\documentclass[a4paper,draft,12pt,reqno]{article}
\usepackage{amsmath,amssymb,amsthm,xspace}
\usepackage{a4wide,amscd,latexsym,amscd}
\usepackage{xy}
\xyoption{all}

\DeclareMathOperator{\Hom}{Hom}
\DeclareMathOperator{\Ext}{Ext}
\newcommand{\bnabla}{\ensuremath{\overline{\nabla}}}
\newcommand{\bDelta}{\ensuremath{\overline{\Delta}}}

\DeclareMathOperator{\rad}{rad}

\DeclareMathOperator{\id}{id}
\DeclareMathOperator{\pd}{pd}

\newcommand{\fd}{\mathcal{F}(\Delta)} 
\newcommand{\fpd}{\mathcal{F}(\overline{\Delta})} 
\newcommand{\fn}{\mathcal{F}(\nabla)} 
\newcommand{\fpn}{\mathcal{F}(\overline{\nabla})} 
\newcommand{\codim}{\operatorname{codim}} 
\newcommand{\tto}{\twoheadrightarrow}

\newcommand{\kk}{\Bbbk}
\renewcommand{\l}{\ensuremath{\lambda}}

\renewcommand{\phi}{\varphi}

\newcommand{\End}{\operatorname{End}}

\newcommand{\ext}{\operatorname{Ext}}
\newcommand{\Add}{\operatorname{Add}}

\newcommand{\mof}{\mathrm{-mod}}
\newcommand{\gl}{\operatorname{gldim}}
\newcommand{\fid}{\operatorname{fdim}}
\newcommand{\fcd}{\operatorname{fcodim}}
\newcommand{\ifid}{\operatorname{ifdim}}

\theoremstyle{plain}
\newtheorem{theorem}{Theorem}
\newtheorem{lemma}{Lemma}
\newtheorem{corollary}{Corollary}

\begin{document} 
\title{On finitistic dimension of stratified algebras}
\author{Volodymyr Mazorchuk}
\date{}
\maketitle

\begin{abstract}
In this survey we discuss the results on the finitistic dimension 
of various stratified algebras. We describe what is already known,
present some recent estimates, and list some open problems. 
\end{abstract}

\section{Introduction and preliminaries}\label{s1}

Let $A$ be a finite-dimensional, associative, and unital algebra over 
an algebraically closed field $\kk$, and $A\mof$ be the category of
finite-dimensional left $A$-modules. Assume that the isomorphism
classes of simple $A$-modules are indexed by $\Lambda=\{1,2,\dots,n\}$
and denote by $L(\l)$, $P(\l)$, $I(\l)$, $\l\in \Lambda$, the
corresponding simple module, its projective cover, and its injective
envelope respectively. Remark that the elements of $\Lambda$ are 
ordered in the natural way. For $\l\in\Lambda$ set $P^{>\l}=
\oplus_{\mu>\l}P(\mu)$ and define the {\em standard} module 
$\Delta(\l)=P(\l)/\mathrm{Trace}_{P^{>\l}}(P(\l))$.
Denote by $\fd$ the full subcategory of $A\mof$, 
which consists of all modules, having
a filtration with subquotients isomorphic to standard modules.
Call the algebra  $A$  {\em strongly standardly stratified} (or an
{\em SSS-algebra})  if ${}_A A\in \fd$. The class of SSS-algebras 
contains the very important subclass of quasi-hereditary 
algebras, and forms a subclass of the class of standardly stratified 
algebras, introduced in \cite{CPS}. SSS-algebras (sometimes also
called just standardly stratified in the literature, which makes
everything somewhat confusing) were intensively studied during 
the last decade, see \cite{AHLU,AHLU2,Ma} and references therein. 
Such algebras arise naturally in Lie theory, see \cite{Ma}. 
In \cite{AHLU} it has been shown that both the
projectively and the injectively defined finitistic dimensions of
such algebras do not exceed $2n-2$. Though this bound is exact for 
certain algebras, in most cases this estimate is very rough. 
For example any hereditary algebra is stratified (even 
quasi-hereditary) with respect to any order on $\Lambda$, see 
\cite[Theorem~1]{DR}, and has  global dimension $1$.

In the present paper we try to approach rather non-symmetric
situations, i.e. the one for which projective and injective 
dimensions can be different. Let $\mathcal{P}^{<\infty}(A)$ and
$\mathcal{I}^{<\infty}(A)$ denote the full subcategories of $A\mof$,
which consists of all modules $M$ having finite projective or
injective dimension respectively. We denote by  $\fid(A)$ the 
{\em projectively defined finitistic dimension} of $A$,
that is the supremum of $\pd(M)$, taken over all $M\in 
\mathcal{P}^{<\infty}(A)$; and by $\ifid(A)$ the {\em injectively 
defined finitistic dimension} of $A$, that is the supremum of
$\id(M)$, taken over all $M\in\mathcal{I}^{<\infty}(A)$.

For $\l\in\Lambda$ define the {\em proper standard} module
$\bDelta(\l)=\Delta(\l)/\mathrm{Trace}_{P(\l)}(\rad\Delta(\l))$. Dually
one defines the {\em costandard} modules $\nabla(\l)$ and the
{\em proper costandard} modules $\bnabla(\l)$, $\l\in\Lambda$. The
categories $\fn$, $\fpd$ and $\fpn$ are defined analogously to $\fd$.
For all modules indexed by $\l\in\Lambda$ the notation without index 
will mean the direct sum over all $\l\in\Lambda$, for example
$L=\oplus_{\l=1}^n L(\l)$ etc. According to \cite{Dl2,La}, an
alternative description of SSS-algebras can be given requiring $I\in
\fpn$. 

Varying the requirements one gets many other classes of stratified
algebras. The ones, which are important for the present paper, are 
{\em properly stratified algebras}, defined in
\cite{Dl} via ${}_A A\in \fd\cap \fpd$, or, alternatively, via
$I\in \fn\cap \fpn$; and  {\em quasi-hereditary algebras}, defined as
those properly stratified algebras, for which $\Delta(\l)=\bDelta(\l)$
for all $\l$, which is equivalent to requiring $\nabla(\l)=\bnabla(\l)$
for all $\l$ (see for example \cite{DR}).

\section{General approach via tilting modules}\label{s2}

\subsection{Tilting modules and finitistic dimension}\label{s2.1}

Let us forget  the stratified structure for a moment. So, let $A$
just be a finite-dimensional, associative, and unital
$\kk$-algebra. Recall, see \cite{Mi}, that an $A$-module $T$ is
called a {\em generalized tilting module} if $T$ has
finite projective dimension, is ext-self-orthogonal, and its additive
closure $\mathrm{Add}(T)$ coresolves ${}_A A$ in a finite number of 
steps. The {\em generalized cotilting modules} are defined dually. 
Looking at the homomorphisms in $D^b(A)$ from $T^{\bullet}[i]$ to the
tilting coresolution of ${}_A A$ one easily derives that $\pd(T)$, in
fact, equals the length of the shortest tilting coresolution of 
${}_A A$. Here for $M\in A\mof$ we denote by $M^{\bullet}$ the 
complex, whose only non-zero component is $M$, concentrated in degree
zero. 

The trivial example of a generalized tilting module is $P$. If 
$\gl(A)<\infty$, then $I$ is a generalized tilting module as well. In
general $I$ need not be a tilting module, since it may 
have infinite projective dimension. However, if $I$ is a
generalized tilting module then, embedding any $M\in
\mathcal{P}^{<\infty}(A)$ into an injective module, and applying
$\Hom_A({}_-,L)$, one derives that $\fid(A)=\pd(I)$. Moreover, in this 
case any $M\in \mathcal{P}^{<\infty}(A)$ can be substituted in
$D^b(A)$ by its finite projective resolution, which then can be
turned into a finite injective complex in $D^b(A)$, since $I$ is a 
tilting module (see for example \cite[Lemma~4]{MO}). This implies 
that any $M\in \mathcal{P}^{<\infty}(A)$ 
has finite injective coresolution, in particular, 
$\mathcal{P}^{<\infty}(A)$ is contravariantly finite in $A\mof$, see
\cite{AR}. 

\subsection{Using self-dual tilting modules}\label{s2.2}

We  have seen that finding non-trivial generalized tilting modules 
in $A\mof$ can give some interesting information about the 
homological behavior of $A\mof$. Especially if such modules are 
self-dual with respect to some contravariant exact equivalence on 
$A\mof$ (usually called a {\em duality}). A duality is called {\em
simple preserving} if it preserve the isomorphism classes of simple
modules. A careful study of the proof of \cite[Theorem~1]{MO} shows 
that what is actually proved there is the following statement:

\begin{theorem}\label{t1}
Let $A$ be a finite-dimensional, associative, and unital
$\kk$-algebra for which $\fid(A)<\infty$. Assume that there exists a 
duality on $A\mof$, and a generalized tilting $A$-module $T$, such
that $Q^{\star}\cong Q$ for every indecomposable $Q\in
\mathrm{Add}(T)$. Then $\fid(A)=2\cdot \pd(T)$. 
\end{theorem}

\begin{proof}
Applying $\star$ to the tilting coresolution of $P$ gives a tilting
resolution of $I$, in particular, $\pd(I)<\infty$. Since
$\fid(A)<\infty$ we can embed any $M\in A\mof$ with $\pd(M)=\fid(A)$
into an injective module, apply $\Hom_A({}_-,L)$, and obtain 
$\pd(I)=\pd(M)=\fid(A)$. Further, $\pd(I)$ is exactly the maximal
degree $l$, for which $\ext_A^l(I,P)$ does not vanish. The latter
can be computed in $D^b(A)$ studying homomorphisms from the shifted
tilting resolution of $I$ to the tilting coresolution of $P$. Under
our assumptions we can apply \cite[Lemma~1]{MO} and the arguments
from \cite[Appendix]{MO}. The statement of the theorem  follows.
\end{proof}

\subsection{Applications to stratified algebras}\label{s2.3}

Assuming $A$ has some sort of stratification makes it in many
cases possible to ensure the assumptions of Theorem~\ref{t1}. 
Indeed, assume that $A$ is an SSS-algebra
having a simple preserving duality (i.e. a duality, which preserves
the isomorphism classes of simple modules). Then $A$ is in fact
properly stratified, the category $\fd\cap\fpn$ equals
$\mathrm{Add}(T)$ for some generalized tilting module $T$
called the {\em characteristic tilting module}, and
the category $\fpd\cap\fn$ equals $\mathrm{Add}(C)$ for some 
generalized tilting module $C$, called the {\em characteristic 
cotilting module}. Moreover, if $T\cong C$, then all
indecomposable direct summands of $T$ are self-dual. The condition
$T\cong C$ is satisfied, for example, for quasi-hereditary algebras. 
Hence we obtain (see \cite[Theorem~1 and Corollary~1]{MO}). 

\begin{corollary}\label{c2}
Let $A$ be an algebra having  a simple preserving duality.
\begin{enumerate}
\item If $A$ is an SSS-algebra and $T\cong C$ then 
$\fid(A)=2\cdot \pd(T)$.
\item If $A$ is quasi-hereditary then  $\gl(A)=2\cdot \pd(T)$.
\end{enumerate}
\end{corollary}

\section{Using tilting and various filtration dimensions}\label{s3}

\subsection{Filtration (co)dimensions}\label{s3.1}

Let $\mathcal{M}$ be a class of $A$-modules and
$\mathcal{F}(\mathcal{M})$ be the full subcategory in $A\mathrm{-mod}$,
which consists of all modules having a filtration with subquotients,
isomorphic to modules from $\mathcal{M}$. For an $A$-module $N$ we
say that $N$ has {\em $\mathcal{M}$-filtration dimension}
(resp. {\em codimension}) $l\in\{0,1,\dots,\infty\}$ if there exists a
resolution (resp. coresolution) of $N$ by modules from 
$\mathcal{F}(\mathcal{M})$ and $l$ is the length of the shortest such
resolution. For properly stratified algebras and SSS-algebras the
following filtration (co)dimensions appear in a natural way:
the {\em Weyl} or {\em standard filtration dimension} $\dim_{\Delta}(N)$ 
for  $\mathcal{M}=\{\Delta(\l),\l\in\Lambda\}$; 
the {\em  proper standard filtration dimension} $\dim_{\bDelta}(N)$ for 
$\mathcal{M}=\{\bDelta(\l),\l\in\Lambda\}$; 
the {\em good} or {\em costandard filtration codimension} 
$\codim_{\nabla}(N)$ for $\mathcal{M}=\{\nabla(\l),\l\in\Lambda\}$; and 
the {\em proper costandard filtration codimension} $\codim_{\bnabla}(N)$
for $\mathcal{M}=\{\bnabla(\l),\l\in\Lambda\}$. If $A$ is an
SSS-algebra, then both $\dim_{\Delta}(N)$ and $\codim_{\bnabla}(N)$
are well-defined for all $N\in A\mof$. In \cite[Lemma~1]{MP} it is
shown that $\dim_{\Delta}(N)=\max\{l|\Ext_A^l(N,\bnabla)\neq 0\}$, 
and $\codim_{\bnabla}(N)=\max\{l|\Ext_A^l(\Delta,N)\neq 0\}$. 
In particular, $\codim_{\bnabla}(N)\leq \pd(\Delta)$ for all $N$, whereas
$\dim_{\Delta}(N)<\infty$ is obviously equivalent to $\pd(M)<\infty$
as $P\in \fd$. We define $\dim_{\Delta}(A)$, $\dim_{\bDelta}(A)$,
$\codim_{\nabla}(A)$, $\codim_{\bnabla}(A)$, $\fid_{\Delta}(A)$ and
$\fcd_{\nabla}(A)$ in the natural way and for an SSS-algebra we obtain
$\codim_{\bnabla}(A)=\pd(\Delta)=\pd(T)$ by \cite[Lemma~1]{MP}.
For properly stratified algebras we dually have 
$\dim_{\bDelta}(A)=\id(\nabla)=\id(C)$. Moreover, by 
\cite[Lemma~2]{MP} we also have 
$\fid_{\Delta}(A)\leq \pd(\nabla)=\pd(C)$ and
$\fcd_{\nabla}(A)\leq \id(\Delta)=\id(T)$.

These filtration (co)dimensions were reinterpreted in 
\cite[Subsection~4.3]{MO} in terms of tilting complexes. Thus we 
have that $\dim_{\Delta}(N)\leq l$ if and only if $N^{\bullet}\in
D^b(A)$ is quasi-isomorphic to a tilting complex 
$\mathcal{T}^{\bullet}$ such that $\mathcal{T}^{i}=0$ for all $i<-l$. 

\subsection{An ``old'' upper bound for $\fid(A)$}\label{s3.2}

The following upper bound  for $\fid(A)$ is stated in  \cite{MP} for 
properly stratified algebras. Here we formulate the result for 
SSS-algebras and present a different proof based on tilting 
resolutions (see also \cite[Corollary~5]{MO}).

\begin{theorem}\label{t3}
Let $A$ be an SSS-algebra. Then $\fid(A)\leq
\fid_{\Delta}(A)+\pd(T)$. 
\end{theorem}

\begin{proof}
If $M\in \mathcal{P}^{<\infty}(A)$, then
$\pd(M)=\max\{l|\Ext_A^l(M,P)\neq 0\}$. We substitute  $M^{\bullet}\in
D^b(A)$ by a quasi-isomorphic tilting complex
$\mathcal{T}^{\bullet}$ satisfying $\mathcal{T}^{i}=0$ for all 
$i<-\dim_{\Delta}(M)$, and we substitute $P$ by its tilting 
coresolution of length $\pd(T)$ (see Subsection~\ref{s2.1}). 
Since for the tilting 
complexes the homomorphisms in $D^b(A)$ can be computed in the 
homotopic category, it is straightforward that $\pd(M)\leq
\dim_{\Delta}(M)+\pd(T)$ and the statement follows.
\end{proof}

If $A$ is properly stratified, as an immediate consequence we have 
$\fid(A)\leq\id(C)+\pd(T)$, which is left-right symmetric and hence 
works for $\ifid(A)$ as well. If $A$ has a duality, everything reduces
to $\fid(A)\leq 2\cdot \pd(T)$. As we have already seen in 
Subsection~\ref{s2.3}, the last bound is exact for quite a wide class 
of quasi-hereditary and stratified algebras, including Schur algebras, 
algebras associated with the BGG-category $\mathcal{O}$ and its
parabolic analogues. 

\subsection{$\fid(A)$ if one can control $\End_A(T)$}\label{s3.3}

Let $A$ be an SSS-algebra. The endomorphism algebra $R=\End_A(T)$ of 
the characteristic tilting module $T$ is called the {\em Ringel dual}
of $A$. The algebra $\End_A(T)^{opp}$ is always an $SSS$-algebra with
respect to the opposite order on $\Lambda$, see 
\cite{AHLU2}. However, $R$ does not need to be properly stratified,
even in the case when  $A$ itself is properly stratified. The algebra
$R$ comes together with the {\em Ringel duality} functor 
$F({}_-)=\Hom_A(T,{}_-):A\mathrm{-mod}\to R\mathrm{-mod}$, which
induces an exact equivalence between the category of $A$-modules 
having a proper costandard filtration and the category of $R$-modules 
having a proper standard filtration.

The Ringel dual $R$ is properly stratified if and only if the module 
$T$ has a filtration with subquotients isomorphic to 
$N(\lambda)=T(\lambda)/\mathrm{Trace}_{T^{<\l}}(T(\l))$, 
where $T^{<\l}=\oplus_{\mu<\lambda}T(\mu)$ (see \cite{FM}). 
In the case when $R$ is properly stratified we denote by $H(\lambda)$,
$\lambda\in\Lambda$, the preimage under $F$ of the indecomposable 
tilting $R$-module corresponding to $\lambda$, and by $H$ the preimage 
under $F$ of the characteristic tilting $R$-module $T^{(R)}$. The 
module $H$ is called the {\em two-step tilting module} for $A$ (since
it is a tilting module for the Ringel dual of $A$). The following 
properties of $H$ were obtained in \cite{FM}:

\begin{theorem}\label{t4}
Assume that $R$ is properly stratified and $H$ is the two-step tilting
module for $A$. Then
\begin{enumerate}
\item $H$ is a generalized tilting module;
\item $\pd(H)=\fid(A)$;
\item $\mathcal{P}^{<\infty}(A)$ coincides with the category of
$A$-modules, which admit a finite coresolution by modules from
$\Add(H)$, in particular, $\mathcal{P}^{<\infty}(A)$ is
contravariantly finite.
\end{enumerate}
\end{theorem}

In particular, the module $H$ is a good test module for $\fid(A)$ and
it completely describes $\mathcal{P}^{<\infty}(A)$ in the homological
sense. It is also shown in \cite{FM} that the existence of $H$ makes it
possible to relate $\fid(A)$ with the projective dimension of the
characteristic tilting module:

\begin{theorem}\label{t5}
Let $A$ be a properly stratified algebra having a simple
preserving  duality. Assume $R$ is properly stratified. Then
\begin{enumerate}
\item $\fid(A)=2\cdot \pd(T^{(R)})$. 
\item $\fid(A)=2\cdot \pd(T)$, in particular, $\pd(T)=\pd(T^{(R)})$,
if $R$ has a simple preserving duality itself.
\end{enumerate}
\end{theorem}

\subsection{A new lower bound for $\fid(A)$}\label{s3.4}

Carefully combining the results of \cite{MO} and \cite{FM} one can
deduce  the following lower bound for the finitistic dimension of 
properly stratified algebras having a simple preserving duality.

\begin{theorem}\label{t8}
Let $A$ be properly stratified with a simple preserving  
duality $\star$. Then we have $\fid(A)\geq 2\cdot \fid_{\Delta}(A)$. 
\end{theorem}

\begin{proof}
We have to produce a module from  
$\mathcal{P}^{<\infty}(A)$ of projective dimension
at least $2\cdot \fid_{\Delta}(A)$. For this it is enough to show that 
any $A$-module $M$, such that $\dim_{\Delta}(M)=\fid_{\Delta}(A)$, 
satisfies $\pd(M)\geq 2\cdot \fid_{\Delta}(A)$. Set
$k=\fid_{\Delta}(A)$. By \cite[Lemma~6]{MO}, $M^{\bullet}$ is
quasi-isomorphic to a finite tilting complex, $\mathcal{T}^{\bullet}$,
satisfying $\mathcal{T}^{i}=0$ for all $i<-k$. Applying $\star$ gives
a finite cotilting  complex $\mathcal{C}^{\bullet}$ satisfying 
$\mathcal{C}^{i}=0$ for all $i>k$. Using \cite[Lemma~11]{FM} one finds
a (possibly infinite) tilting complex $\mathcal{Q}^{\bullet}$, which
is quasi-isomorphic to $\mathcal{C}^{\bullet}$, and which
satisfies $\mathcal{Q}^{i}=0$ for all $i>k$. Moreover, using 
\cite[Lemma~12]{FM} one can also guarantee that $\mathcal{T}^{-k}$ is
non-trivial and is a direct summand of $\mathcal{Q}^{k}$. Using the
arguments as in  \cite[Section~3]{MO} one shows that
there is a non-zero morphism  from $\mathcal{T}^{\bullet}[-2k]$ to 
$\mathcal{Q}^{\bullet}$, implying $\pd(M)\geq 2k$.
\end{proof}

It is interesting to compare the bound, given in Theorem~\ref{t8},
with the results, described in Subsection~\ref{s3.3}. For this we will
need the following lemma:

\begin{lemma}\label{l9}
Let $A$ be an SSS-algebra and $M\in\fpn$ such that $\pd(M)<\infty$. 
Then $\dim_{\Delta}(M)=\pd(F(M))$.
\end{lemma}

\begin{proof}
Taking the minimal projective resolution $\mathcal{P}^{\bullet}$ of
$M$ and applying \cite[Lemma~4.1]{MO} we obtain a finite tilting
complex $\mathcal{T}^{\bullet}$, which is quasi-isomorphic to 
$M^{\bullet}\in D^b(A)$. Using the arguments from the proof of 
\cite[Lemma~5]{MO} one even shows that $\mathcal{T}^{\bullet}$ is 
quasi-isomorphic to  a finite minimal (in the sense of \cite{MO})
tilting complex $\mathcal{Q}^{\bullet}$ 
satisfying $\mathcal{Q}^{i}=0$, $i>0$. In other words, the module $M$
admits a finite tilting resolution. Applying $F$ gives a projective
resolution of $F(M)$ and we see that the length of the minimal tilting 
resolution of $M$ is exactly $\pd(F(M))$. From \cite[Lemma~6]{MO}
it also follows that the length of the minimal tilting resolution of
$M$ equals $\dim_{\Delta}(M)$, completing the proof.
\end{proof}

\begin{corollary}\label{c10}
Let $A$ be properly stratified and assume that $R$ is also
properly stratified. Then $\pd(T^{(R)})=\fid_{\Delta}(A)$.
\end{corollary}

\begin{proof}
By Lemma~\ref{l9} we have $\pd(T^{(R)})=\dim_{\Delta}(H)$. Further,
let $M\in \mathcal{P}^{<\infty}(A)$ be such that 
$l=\dim_{\Delta}(M)=\fid_{\Delta}(A)$. By Theorem~\ref{t4}, we have a 
short exact sequence $M\hookrightarrow H_1\tto K$, where
$H_1\in\Add(H)$ and $K\in 
\mathcal{P}^{<\infty}(A)$. In particular, $\dim_{\Delta}(H_1)$ and 
$\dim_{\Delta}(K)$ do not exceed $\dim_{\Delta}(M)$. Applying
$\Hom_A({}_-,\bnabla)$ we obtain that $\Ext_A^{l}(H_1,\bnabla)$ surjects
onto $\Ext_A^{l}(M,\bnabla)\neq 0$ and hence $\dim_{\Delta}(H_1)=l$ by
\cite[Lemma~1]{MP}. This implies $\dim_{\Delta}(H)=l$ and completes
the proof. 
\end{proof}

An immediate corollary of Theorem~\ref{t5} and Corollary~\ref{c10} is:

\begin{corollary}\label{c11}
Let $A$ be properly stratified having a simple preserving duality. 
Assume that $R$ is also properly stratified. Then 
$\fid(A)=2\cdot \fid_{\Delta}(A)$.
\end{corollary}

\section{A counterexample}\label{s4}

In \cite[Conjecture~1]{MP} it was conjectured that the finitistic
dimension of a properly stratified algebra having a simple preserving 
duality always equals twice the projective dimension of the 
characteristic tilting module. As we saw above this is true under 
assumptions that $R$ is properly stratified and has a simple
preserving duality, which includes, in particular, the cases
of quasi-hereditary algebras, and properly stratified algebras whose
tilting modules are also cotilting. Unfortunately, in the full
generality the statement of the conjecture is wrong. As a
counter example one can consider the following algebra (the first 
counter example was constructed by the author, computed by Birge 
Huisgen-Zimmermann, and simplified by Steffen K{\"o}nig).

Let $A$ be the path algebra of the quiver
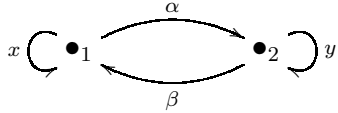
\begin{displaymath}
\xymatrix{
\bullet_{1} \ar@/^1pc/[rr]^{\alpha}\ar@(ul,dl)[]_{x} 
& &  
\bullet_{2}
\ar@/^1pc/[ll]^{\beta} \ar@(ur,dr)[]^{y}   
}
\end{displaymath}
modulo the relations $\alpha\beta=x^2=y^2=x\beta=\alpha x=0$.
The map $\alpha\mapsto \beta$, $\beta\mapsto\alpha$ extends to an  
anti-involution on $A$ and hence gives rise to a duality on
$A\mathrm{-mod}$.  

The radical filtrations of the projective, standard, and proper 
standard modules look as follows:
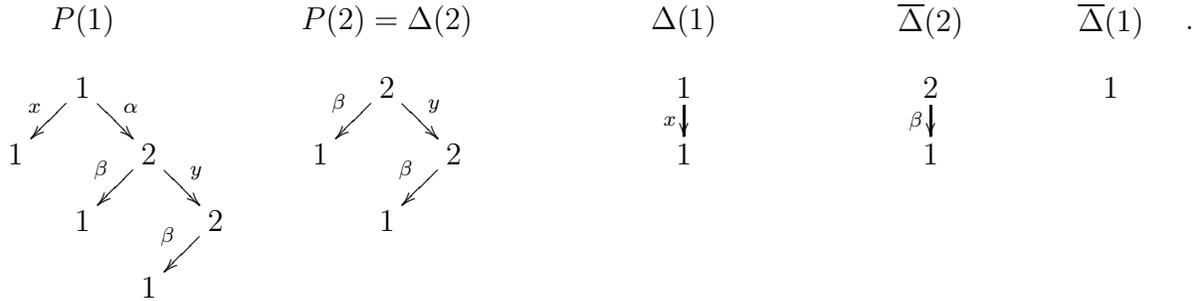
\begin{displaymath}
\xymatrix@!=0.1pc{
  & P(1)  &  &  \\
& 1\ar@{->}[rd]^{\alpha}\ar@{->}[ld]_{x} &   &  \\
 1 &  & 2\ar@{->}[rd]^{y}\ar@{->}[ld]_{\beta}   &  \\
 & 1 &   & 2 \ar@{->}[ld]_{\beta} \\
 &   & 1 &  
}
\qquad
\xymatrix@!=0.1pc{
  & P(2)=\Delta(2)  &  \\
& 2\ar@{->}[rd]^{y}\ar@{->}[ld]_{\beta}   &  \\
1 &   & 2 \ar@{->}[ld]_{\beta} \\
  & 1 &   \\
 &   &    
}
\qquad
\xymatrix@!=0.1pc{
&   & \Delta(1)  &  \\
&   & 1\ar@{->}[d]_{x} &     \\
&    & 1 &     \\
&    &   &     \\
&    &   &   & \\
&    &   &   
}
\quad
\xymatrix@!=0.1pc{
  & \bDelta(2)  &  \\
 & 2\ar@{->}[d]_{\beta} &     \\
  & 1 &     \\
  &   &     \\
  &   &    \\
  &   &   
}\quad
\xymatrix@!=0.1pc{
  & \bDelta(1)  &  \\
 & 1 &     \\
  &  &     \\
  &   &     \\
  &   &    \\
  &   &   
}.
\end{displaymath}
It follows that 
$A$ is properly stratified. Since $A$ has a duality, the socle
filtrations of injective, costandard and proper costandard modules
are duals of the radical filtrations of the corresponding projective,
standard and proper standard modules above. The indecomposable
tilting $A$-modules have the following radical filtration:
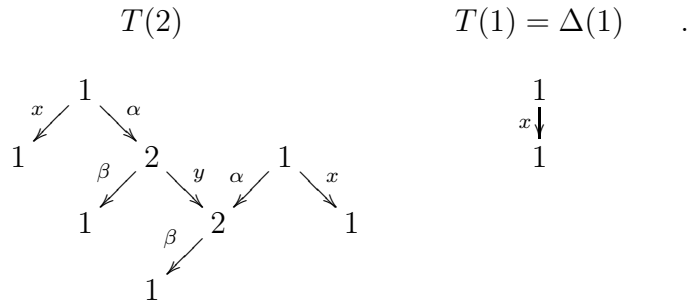
\begin{displaymath}
\xymatrix@!=0.1pc{
 &   & T(2)  &  &  & \\
 & 1\ar@{->}[ld]_{x}\ar@{->}[rd]^{\alpha}  &   &  &  & \\
1  &   & 2\ar@{->}[ld]_{\beta}\ar@{->}[rd]^{y}  &  & 
1\ar@{->}[ld]_{\alpha}\ar@{->}[rd]^{x} & \\
 &  1 &   & 2\ar@{->}[ld]_{\beta}  &  & 1\\
 &   &  1 &  &  & \\
}\quad
\xymatrix@!=0.1pc{
&   & T(1)=\Delta(1)  &  \\
&   & 1\ar@{->}[d]_{x} &     \\
&    & 1 &     \\
&    &   &     \\
&    &   &   & \\
&    &   &   
}
.
\end{displaymath}
Neither $T(2)$ nor $\Delta(1)$ are projective, which
implies that $\fid(A)\geq 1$. Further it is easy to see that there are
the following minimal projective resolutions of tilting modules: 
$0\to P(2)\to P(1)\to T(1)\to 0$ and $0\to P(2)\to P(1)\oplus
P(1)\to T(2)\to 0$, and hence $\pd(T)=1$. It is also easy to see that 
any injection between tilting modules is an isomorphism. This and 
\cite[Lemma~6]{MO} implies $\fid_{\Delta}(A)=0$. Hence, by
Theorem~\ref{t3} we obtain $\fid(A)\leq \fid_{\Delta}(A)+\pd(T)=1$ and
thus $\fid(A)=1$. In particular, in this example we have 
$2\cdot\fid_{\Delta}(A)=0 <\fid(A)=1< 2\cdot \pd(T)=2$.

The example together with Theorem~\ref{t3} motivates to make the 
following correction to \cite[Conjecture~1]{MP}:
\vspace{0.5cm}

\noindent
{\bf Corrected conjecture.} {\em Let $A$ be a properly stratified algebra
with a simple preserving duality. Then
$\fid(A)=\fid_{\Delta}(A)+\pd(T)$.} 
\vspace{0.5cm}

Remark that for algebras having a simple preserving duality 
we always have $\fid_{\Delta}(A)\leq \pd(T)=\codim_{\bnabla}(A)$, 
see \cite{MP}.

\section{A bound for $\ifid(A)$ in the case of an SSS-algebra}\label{s5}

Up to this point all the results mentioned were about
the projectively defined finitistic dimension. A natural question
is: what can one say about the injectively defined version? In the
case of an algebra having a (simple preserving) duality the answer 
is very easy: the
injectively and the projectively defined finitistic dimensions
obviously coincide. But what can be said in the general case? This
question is still more or less open, see Section~\ref{s6}. Here we
just present an easy upper bound for the case, when we have enough
information about the Ringel dual of the algebra.

\begin{theorem}\label{t12}
Let $A$ be and SSS-algebra and assume that $R$ is properly stratified.
Then $\ifid(A)\leq \fid(A)\leq \fid_{\Delta}(A)+\pd(T)$. 
\end{theorem}

\begin{proof}
Because of Theorem~\ref{t3} it is enough to prove that
$\ifid(A)\leq \fid(A)$. To do this we will show that for any
$M$ with $\id(M)=k<\infty$ we have $\Ext_A^k(H,M)\neq 0$.
From the definition of $H(\lambda)$ it follows that there is an
injection $\bnabla(\lambda)\hookrightarrow H(\lambda)$, and hence
there is an injection $L\hookrightarrow H$ with cokernel, $K$
say. Applying $\Hom_A({}_-,M)$ to the short exact sequence
$L\hookrightarrow H\tto K$ and using $\id(M)=k$ we get a surjection of 
$\Ext_A^k(H,M)$ onto $\Ext_A^k(L,M)\neq 0$. This completes the proof.
\end{proof}

\section{Some comments and questions}\label{s6}

Summarizing the results of the paper  we can say the following: if
$A$ is a properly stratified algebra having a simple preserving
duality, then we have the following bounds for $\fid(A)$:
\begin{equation}\label{eq1}
2\cdot\fid_{\Delta}(A)\leq \fid(A)\leq \fid_{\Delta}(A)+\pd(T)\leq
2\cdot \pd(T).
\end{equation}
Most of the components of \eqref{eq1} have equivalent reformulations
in other homological terms for $A$ or for the Ringel dual $R$ of $A$.
In many cases, for example for quasi-hereditary algebras, we know that 
all inequalities in \eqref{eq1} are in fact equalities. 
We also know that the first and the third inequalities can be strict. 
This gives rise to the following question: 

\begin{enumerate}
\item Let $A$ be an SSS-algebra. How different can $\pd(T)$ and 
$\fid_{\Delta}(A)$ be? The same for properly stratified algebras
and for properly stratified algebras with duality.
\item Describe the class of SSS-algebras with properly stratified
Ringel duals, satisfying $\pd(T)=\pd(T^{(R)})$ (remark that the last
condition immediately makes all inequalities of \eqref{eq1} into
equalities). The same for properly stratified algebras and for 
properly stratified algebras with duality.
\end{enumerate}

We saw  that the module $H$, which appears in the case when $R$ is
properly stratified, can be used as a test module for $\fid(A)$. It
was shown in \cite{FM} that $\End_A(H)^{opp}$ is always an
SSS-algebra. Hence, very natural questions are:

\begin{enumerate}
\setcounter{enumi}{2}
\item Find, in terms of $A\mathrm{-mod}$, necessary and sufficient 
conditions for $\End_A(H)$ to be properly stratified.
\item Is there any relation between $\fid(A)$ and $\fid(\End_A(H))$? 
\end{enumerate}

As we have already mentioned, much less is know about the injectively 
defined finitistic dimension, so even the following very general
question is not answered: 

\begin{enumerate}
\setcounter{enumi}{4}
\item Let $A$ be an SSS-algebra or a properly stratified algebra. 
Can one use tilting modules to estimate or compute $\ifid(A)$?
\end{enumerate}

\begin{flushleft}
\bf Acknowledgments
\end{flushleft}

The author acknowledges the supports of The Swedish Research Council, 
The Royal Swedish Academy of Sciences, and The Swedish Foundation for
International Cooperation in Research and Higher Education (STINT).
The author thanks Anders Frisk  and Catharina Stroppel for their 
comments on the paper.

\vspace{1cm}

\noindent
Volodymyr Mazorchuk, Department of Mathematics, Uppsala University,
Box 480, 751 06, Uppsala, SWEDEN, 
e-mail: {\tt mazor\symbol{64}math.uu.se},
web: {``http://www.math.uu.se/$\tilde{\hspace{1mm}}$mazor/''}.
\vspace{0.5cm}

\end{document}